\documentclass{amsart}

\usepackage{amsthm}
\usepackage{amsmath}
\usepackage{amsfonts}
\usepackage{fullpage}
\usepackage{pinlabel}
\usepackage{pictex}
\usepackage{graphicx}
\usepackage{etex}
\usepackage{amscd}

\usepackage{hyperref}  
\hypersetup{%
  bookmarksnumbered=true,%
  bookmarks=true,%
  colorlinks=true,%
  linkcolor=blue,%
  citecolor=blue,%
  filecolor=blue,%
  menucolor=blue,%
  pagecolor=blue,%
  urlcolor=blue,%
  pdfnewwindow=true,%
  pdfstartview=FitBH}
\def\arxiv#1{\relax\ifhmode\unskip\quad\fi
    \href{http://arxiv.org/abs/#1}%
{\tt arXiv:\penalty -100\unskip#1}}    

\def\MR#1{\relax\ifhmode\unskip\quad\fi
    \href{http://www.ams.org/mathscinet-getitem?mr=#1}{MR#1}}
\def\xox#1{\csname xx#1\endcsname} 
\def\xxZBL#1{}
\def\xxJFM#1{}
\def\@url#1{{\tt\def~{\lower3.5pt\hbox{\char'176}}\def\_{\char'137}#1}}
\let\fullref\autoref
\def\makeautorefname#1#2{\expandafter\def\csname#1autorefname\endcsname{#2}}
\makeautorefname{figure}{Figure}%
\makeautorefname{section}{Section}%
\makeautorefname{subsection}{Section}%
\makeautorefname{theorem}{Theorem}%
\makeautorefname{thm}{Theorem}%
\makeautorefname{corollary}{Corollary}%
\makeautorefname{lemma}{Lemma}%
\makeautorefname{prop}{Proposition}%
\makeautorefname{definition}{Definition}%
\makeautorefname{rem}{Remark}%
\makeautorefname{fact}{Fact}%

\title[Examples of Non-Rigid CAT(0) Groups from the Category of Knot Groups]{Examples of Non-Rigid CAT(0) Groups from the Category of Knot Groups}

\author{Christopher Mooney}
\address{Department of Mathematical Sciences\\
University of Wisconsin-Milwaukee\\
PO Box 413\\
Milwaukee, Wisconsin 53201-0413}
\email{cpmooney@uwm.edu}
\date{6 August 2008}

\newtheorem{prop}{Proposition}[section]

\newtheorem{fact}{Fact}[section]
\makeatletter
\let\c@fact=\c@prop
\makeatother

\newtheorem{lemma}{Lemma}[section]
\makeatletter
\let\c@lemma=\c@prop
\makeatother

\newtheorem{corollary}{Corollary}[section]
\makeatletter
\let\c@corollary=\c@prop
\makeatother

\newtheorem{rem}{Remark}[section]
\makeatletter
\let\c@rem=\c@prop
\makeatother

\newtheorem{definition}{Definition}[section]
\makeatletter
\let\c@definition=\c@prop
\makeatother

\newtheorem*{ThmCK}{Theorem CK}
\newtheorem*{ThmW}{Theorem W}

\newtheorem*{Thm1}{Theorem 1}
\newtheorem*{ThmA}{Theorem A}
\newtheorem*{ThmB}{Theorem B}
\newtheorem*{ThmC}{Theorem C}
\newtheorem*{ThmD}{Theorem D}
\newtheorem*{ThmE}{Theorem E}

\newtheorem*{ThmA'}{Theorem A$'$}
\newtheorem*{ThmB'}{Theorem B$'$}
\newtheorem*{ThmC'}{Theorem C$'$}
\newtheorem*{ThmD'}{Theorem D$'$}
\newtheorem*{ThmE'}{Theorem E$'$}

\newtheorem*{claim}{Claim}

\DeclareMathOperator{\im}{im}

\DeclareMathOperator{\Itin}{Itin}

\DeclareMathOperator{\Min}{Min}

\newcommand{\incl}{\hookrightarrow}

\newcommand{\N}{\mathcal{N}}
\newcommand{\C}{\mathcal{C}}

\newcommand{\Z}{\mathbb{Z}}
\newcommand{\R}{\mathbb{R}}

\begin{document}

\begin{abstract}
C Croke and B Kleiner have constructed an example of a CAT(0) group with more than one visual boundary.  J Wilson has proven that this same group has uncountably many distinct boundaries.  In this article we prove that the knot group of any connected sum of two non-trivial torus knots also has uncountably many distinct CAT(0) boundaries.
\end{abstract}

\maketitle

\section{Introduction}
\label{sec:intro}
The CAT(0) condition is a geometric notion of nonpositive curvature similar to the
definition of Gromov $\delta$--hyperbolicity.  A proper geodesic space $X$ is called CAT(0)
if it has the property that geodesic triangles in $X$ are ``no fatter'' than geodesic
triangles in Euclidean space (a precise definition is given by M\,R Bridson and A Haefliger in
\cite[Chapter II.1]{BH}).
The \textit{visual} or \textit{ideal boundary} of $X$, denoted \(\partial X\),
is the collection of endpoints of geodesic rays emanating from a chosen basepoint
endowed with the cone topology.
It is well-known that \(\partial X\) is well-defined and independent of choice of basepoint
and that \(X\cup\partial X\) is a $Z$--set compactification for $X$.
A group $G$ is called CAT(0) if it acts geometrically
(i.e. properly discontinuously and cocompactly by isometries) on some CAT(0) space $X$.
In this setup we call $X$ a CAT(0) $G$--space
and \(\partial X\) a CAT(0) boundary of $G$.
We say that a CAT(0) group $G$ is \textit{rigid} if it has only one
topologically distinct boundary.\\

It is well-known that if $G$ is negatively curved (acts geometrically on a Gromov
$\delta$--hyperbolic space) or if $G$ is free abelian then $G$ is rigid.  Apart from
this little is known concerning rigidity of groups.  P\,L Bowers and K Ruane
showed that if $G$ splits as the product of a negatively curved group with
a free abelian group then $G$ is rigid \cite{BR}.
Ruane proved later in \cite{Ru} that if $G$ splits as a product of two negatively curved groups then $G$
is rigid.  T Hosaka has extended this work to show that in fact it suffices to know
that $G$ splits as a product of rigid groups \cite{Ho}.
Another condition which guarantees
rigidity is knowing that $G$ acts on a CAT(0) space with isolated flats which was proven
by C Hruska in \cite{Hr}.\\

Not all CAT(0) groups are rigid, however: C Croke and B Kleiner constructed in \cite{CK}
an example of a non-rigid CAT(0) group $G$.  Specifically, they showed that $G$ acts on two different
CAT(0) spaces whose boundaries admit no homeomorphism.
J Wilson proved in \cite{Wi} that this same group has \textit{uncountably many}
boundaries.\\

In this article we exhibit an infinite family of non-rigid knot groups.
It is a Corollary of Thurston's hyperbolization theorem \cite{Th}
for Haken 3--manifolds that every knot is either a torus knot, a hyperbolic knot
or a satellite knot.  It follows from Hruska's result \cite{Hr} that hyperbolic knot groups are rigid.
Furthermore, using a result of T Bedenikovic, A Delgado and M Timm \cite{BDT}
and the Bowers--Ruane result from \cite{BR} we can prove
that torus knot groups are rigid (see \fullref{prop:rigidtorusknots}).
The following theorem gives us an infinite family of non-rigid satellite knots.\\

\begin{Thm1}
The knot group $G$ of any connected sum of two non-trivial torus knots has
uncountably many CAT(0) boundaries.
\end{Thm1}

Specifically, we will prove that given any such knot group $G$,
there is a natural construction of a family of CAT(0) $G$--spaces which is analogous
to the construction used by Croke and Kleiner in \cite{CK}.
Even though each space here will have a similar but significantly different structure from
the spaces constructed in \cite{CK} (see \fullref{sec:connectedsumsoftorusknots}), we will
show that on the level of boundaries they have the same basic properties.
Interestingly enough, it turns out that the proof given in \cite{CK} will not
work in this situation.  In order to get any results we will require the work of
Wilson \cite{Wi}.
This is discussed in more detail at the end of \fullref{sec:origCK}.\\

As a final comment on the statement of Theorem 1,
there is a stronger notion of rigidity than the definition we use here.
Sometimes a CAT(0) group is said to be \textit{rigid} if every $G$--equivariant quasi-isometry
between two CAT(0) spaces extends to a homeomorphism of the boundaries.
For us such a group will be called \textit{strongly rigid}.
Negatively curved groups are strongly rigid, for instance.  The fact that these two notions
of rigidity are distinct
is due to Bowers and Ruane who exhibit in \cite{BR} an example of a group which is rigid
(that is, \textit{weakly} rigid) but not strongly rigid.\\

In \cite{CK2} Croke and Kleiner found necessary and sufficient conditions for determining
when the fundamental group of a 3--dimensional graph manifold is strongly rigid.
Since the groups we are considering fall under this category,
our result is slightly stronger than theirs for this particular class of groups.
We prove that the knot group of any connected sum of two non-trivial torus knots
is not even weakly rigid.\\

\subsection{Acknowledgements}
The work contained in this paper is published as one part of the author's Ph.D. thesis written
under the direction of Craig Guilbault at the University of Wisconsin-Milwaukee.
The author would like to thank the referee for their helpful suggestions along with
Ric Ancel, Chris Hruska, Boris Okun and Tim Schroeder.

\section{Croke and Kleiner's Original Construction}
\label{sec:origCK}
Before diving into the proof of Theorem 1, we quickly sketch the proof of the main theorem of \cite{CK}.
Let $G=G_{CK}$ be the group given by the presentation:
\[
	\left<a,b,c,d|ab=ba,bc=cb,cd=dc\right>
\]
Croke and Kleiner construct CAT(0) $G$--spaces $X$ such that each $X$ is covered by a collection of closed convex subspaces
called \textit{blocks}.  The visual
boundary $\partial B$ of every block $B$ is the suspension of a Cantor set.  The suspension points
are called \textit{poles}.
If two blocks $B_0$ and $B_1$ intersect, then $B_0$ is said to \textit{neighbor} $B_1$ and their intersection is a
Euclidean plane called a \textit{wall}.
They then prove five statements for each $X$\\

\begin{ThmA}{\rm \cite[Section 1.4]{CK}}\quad
The nerve $N$ of the collection of blocks is a tree.
\end{ThmA}

\begin{ThmB}{\rm\cite[Lemma 3]{CK}}\quad
Let \(B_0\) and \(B_1\) be blocks and $D$ be the distance between the corresponding vertices in
$N$.  Then:\\
\begin{enumerate}
\item If \(D=1\), then \(\partial B_0\cap\partial B_1=\partial W\) where $W$ is the wall \(B_0\cap B_1\).\\
\item If \(D=2\), then \(\partial B_0\cap\partial B_1\) is the set of poles of \(B_{\frac12}\) where
\(B_{\frac12}\) intersects \(B_0\) and \(B_1\).\\
\item If \(D>2\), then \(\partial B_0\cap\partial B_1=\emptyset\).\\
\end{enumerate}
\end{ThmB}

A \textit{local path component} of a point in a space is a path component of an open neighborhood of that point.\\

\begin{ThmC}{\rm \cite[Lemma 4]{CK}}\quad
Let $B$ be a block and \(\zeta\in\partial B\) not be a pole of any neighboring block.  Then $\zeta$ has a local
path component which stays in \(\partial B\).
\end{ThmC}

\begin{ThmD}{\rm \cite[Corollary 8]{CK}}\quad
The union of block boundaries in $\partial X$ is the unique dense \textit{safe path} component
of \(\partial X\).
\end{ThmD}

The definition of \textit{safe path} will be given in
\fullref{sec:KnotGroupBoundary}.
For now it suffices to understand that
Theorem D gives a way to topologically distinguish the union of block boundaries in $\partial X$.  With these
thereoms in hand it is not hard to prove that given two constructions $X_1$ and $X_2$, any homeomorphism
\(\partial X_1\to\partial X_2\) takes poles to poles, block boundaries to block boundaries and wall boundaries
to wall boundaries.  The last piece of the puzzle is Theorem E.  Given \(0<\theta\le\pi/2\), we can construct
$X_{\theta}$ in such a way that the minimum Tits distance between poles is $\theta$.  For a block $B$, we denote by $\Pi B$
the set of poles of neighboring blocks.\\

\begin{ThmE}{\rm \cite[Lemma 9]{CK} (also \cite[Proposition 2.2]{Wi})}\quad
For a block $B$, the union of boundaries of walls of $B$ is dense in \(\partial B\) and
\(\overline{\Pi B}\) is precisely the set of points of \(\partial B\) which are a Tits distance
of $\theta$ from a pole of $B$.
\end{ThmE}

With these five theorems in hand we get the main result of \cite{CK}:\\

\begin{ThmCK}
Let $B$ be a block and $L$ be a suspension arc of \(\partial B\).  Then \(\bigl|L\cap\overline{\Pi B}\bigr|=1\) iff
\(\theta=\pi/2\).  Therefore $G_{CK}$ has at least two distinct boundaries.
\end{ThmCK}

In \cite{Wi} Wilson uses these five theorems to prove a stronger result:\\

\begin{ThmW}
If \(\theta_1\neq\theta_2\), then \(\partial X_{\theta_1}\not\approx\partial X_{\theta_2}\).
Therefore $G_{CK}$ has uncountably many distinct boundaries.
\end{ThmW}

In this article we consider the knot group $G=G_K$ of any connected sum $K$ of torus knots.
We produce for $G$ an analogous family of CAT(0) spaces which have a similar structure to those constructed in \cite{CK}.
Specifically, we have \textit{blocks}, \textit{walls} and \textit{poles} for these spaces as well,
and for each \(0<\theta<\pi/2\) we can construct \(X_\theta\) such that
the minimum Tits distance between two poles is $\theta$.  This done, we show that we have the
appropriate analogues to Theorems A--E.\\

Now if we had \(X_{\pi/2}\), then Theorems A--E would be enough to guarantee that $G$
has at least two boundaries.  Thus we would not need the arguments found in \cite{Wi}
to prove that $G$ is not rigid.
However, as we will see in \fullref{prop:nopiover2} there is no ``natural'' construction which
will yield \(X_{\pi/2}\).  Therefore in order to prove that $G$ is not rigid we really need
to apply the work of \cite{Wi}.\\

\section{Block Structures on CAT(0) Spaces}
\label{sec:NervesandItineraries}
We begin by observing that the work in Sections 1.4--5 of \cite{CK} does not depend on the specific
construction used in in \cite{CK}.  The same observations apply if we replace their definition of a
\textit{block} with the following one.\\

\begin{definition}
\label{defn:blockstructure}
Let $X$ be a CAT(0) space and \(\mathcal{B}\) be a collection of closed convex subspaces covering $X$.
We call $\mathcal B$ a \textit{block structure on $X$} and its elements \textit{blocks} if
\(\mathcal{B}\) satisfies the following three properties:\\
\begin{enumerate}
\item Every block intersects at least two other blocks.\\
\item Every block has a $(+)$ or $(-)$ parity such that two blocks intersect only if they have opposite parity.\\
\item There is an \(\epsilon>0\) such that two blocks intersect iff their $\epsilon$--neighborhoods intersect.\\
\end{enumerate}
\end{definition}

The \textit{nerve} of a collection \(\mathcal{C}\) of sets is
the (abstract) simplicial complex with vertex set \(\{v_B|B\in\mathcal{C}\}\)
such that a simplex \(\{v_{B_1},...,v_{B_n}\}\) is included whenever
\(\bigcap_{i=1}^nB_i\neq\emptyset\).
In exactly the same way as in \cite{CK} the nerve $\N$ of the collection of blocks
is a tree, and we can define the itinerary of a geodesic.
A geodesic $\alpha$ is said to \textit{enter} a block if it passes through a point which is not in any other block.
The \textit{itinerary of $\alpha$} is defined to be the list
\([B_1,B_2,...]\) where $B_i$ is the $i^{\textrm{th}}$ block that $\alpha$ enters.
This list is denoted by \(\Itin\alpha\).
The following lemma follows in the same way as \cite[Lemma 2]{CK},
which simply uses the fact that a block $B$ is convex and that
its topological frontier is covered by the collection of blocks corresponding to the link
in $\N$ of the vertex $v_B$.\\

\begin{lemma}
\label{le:itingeo}
If \(\Itin\alpha=[B_1,B_2,...]\), then \([v_{B_1},v_{B_2},...]\) is a geodesic in $\N$.\\
\end{lemma}

We may also talk about the \textit{itinerary between two blocks}.
If \([v_{B_1},...,v_{B_n}]\) is the geodesic edge path in $\N$ connecting
two vertices \(v_{B'_0}\) and \(v_{B'_1}\), then we call \([B_1,...,B_n]\) the itinerary between
\(B'_0\) and \(B'_1\) and write:
\[
	\Itin[B'_0,B'_1]=[B_1,...,B_n]
\]
The two notions of itineraries are related as follows:
The itinerary of a geodesic segment $\alpha$ is the shortest itinerary \(\Itin[B'_0,B'_1]\) for which $\alpha$ begins in
$B'_0$ and ends in $B'_1$.  Note also that the same observations which gave us \fullref{le:itingeo}
also provide this next lemma.\\

\begin{lemma}
\label{le:itins}
Let $B'_0$ and $B'_1$ be blocks, write \(\Itin[B'_0,B'_1]=[B_1,...,B_n]\), and let $\alpha$ be a geodesic
beginning in $B'_0$ and ending in $B'_1$.
Then:
\begin{enumerate}
\item $\alpha$ enters $B_k$ for every \(1<k<n\).\\
\item $\alpha$ passes through \(B_k\cap B_{k+1}\) for every \(1\le k<n\).\\
\item \(\bigcup_{k=1}^nB_k\) is convex.\\
\end{enumerate}
\end{lemma}

We call a geodesic ray \textit{rational} if its itinerary is finite and
\textit{irrational} if its itinerary is infinite.  A point of \(\partial X\)
is called \textit{irrational} if it is the endpoint of an irrational geodesic
ray; otherwise we call it \textit{rational}.  We denote the set of rational points
of \(\partial X\) by $RX$ and the set of irrational points by $IX$.\\

\begin{lemma}
\label{le:irrationaldistance}
Let $\alpha$ be an irrational geodesic ray.
Then for any block \(B_0\)
\[
	\lim_{t\to\infty}d\bigl(\alpha(t),B_0\bigr)=\infty.
\]
\end{lemma}

\begin{proof}
Write \(\Itin\alpha=[B_1,B_2,...]\).  Since $\N$ is a tree we can find \(M>1\)
such that for every \(m\ge M\), \(\Itin[B_0,B_m]\ni B_M\).
For $m\ge M$ choose a time $t_m$ such that \(\alpha(t_m)\in B_m\).
Then:
\begin{align*}
	\lim_{t\to\infty}
	d\bigl(\alpha(t),B_0\bigr)
&\ge
	\lim_{t\to\infty}
	d\bigl(\alpha(t),B_M\bigr) \\
&=
	\lim_{m\to\infty}
	d\bigl(\alpha(t_m),B_M\bigr) \\
&\ge
	\lim_{m\to\infty}
	d\bigl(B_m,B_M\bigr)
\end{align*}
Hence it suffices to prove the following.\\

\begin{claim}
Let $\epsilon$ be given as in condition (3) of \fullref{defn:blockstructure}.
Then whenever \(d(v_{B},v_{B'})\ge 2k\), we have \(d(B,B')\ge 2k\epsilon\).
\end{claim}

Note that whenever \(d(v_{B},v_{B'})=2\) then we have \(d(B,B')\ge 2\epsilon\)
because the $\epsilon$--neighborhoods of $B$ and $B'$ do not overlap.
Assume \(\Itin[B,B']=[B_0,B_1,...,B_n]\)
where \(n\ge 2k\).  Then for any \(x\in B\) and \(x'\in B'\)
the geodesic \([x,x']\) passes through
\(B_{2i}\) for \(0\le i\le k\) at some point $z_i$.
So
\[
	d(x,x')
=
	\sum_{i=0}^{k-1}d(z_i,z_{i+1})
\ge
	2k\epsilon.
\]
\end{proof}

\begin{corollary}
\label{co:RXIX}
\rule{0cm}{0cm}
\begin{enumerate}
\item \(RX\) is the union of block boundaries in \(\partial X\) and $IX$ is its complement.\\
\item If \(\zeta\in IX\), then every geodesic ray going out to $\zeta$ is irrational.\\
\item If \(\zeta\in IX\) and $\alpha$ and $\beta$ are geodesic rays going out to $\zeta$,
then the itineraries of $\alpha$ and $\beta$ eventually coincide.\\
\end{enumerate}
\end{corollary}

A geodesic space is said to have the \textit{geodesic extension property} if every geodesic segment
can be extended to a geodesic line.  As is true with the original Croke--Kleiner construction, the
blocks we construct will satisfy the geodesic extension property.\\

\begin{lemma}
\label{le:denseRX}
If blocks have the geodesic extension property, then $RX$ is dense.
\end{lemma}

\begin{proof}
Let $\alpha$ be an irrational geodesic ray and write \(\Itin\alpha=[B_1,B_2,...]\).
For each \(n\ge 1\) let $t_n$ be a time at which \(\alpha(t_n)\in B_n\).
Then every ray \(\alpha|_{[0,t_n]}\) can be extended to a geodesic ray \(\alpha_n\)
which does not leave the block \(B_n\).  Then \(\alpha_n\to\alpha\).
\end{proof}

We end this section with a definition which will simplify the proof of Theorem D$'$ later.
Given a space $Y$ we call a surjective map \(\phi: IX\to Y\) an \textit{irrational map} if it satisfies the
property that \(\phi(a)=\phi(b)\) iff whenever $\alpha$ and $\beta$ are geodesic rays
going out to \(a\) and \(b\) respectively then \(\Itin\alpha\) and \(\Itin\beta\)
eventually coincide.
The obvious candidate for such a map is the function \(\phi: IX\to\partial\N\) which takes
$a$ to the boundary point in $\partial\N$ determined by the itinerary of a ray going out to $a$.
This function is well-defined by \fullref{co:RXIX}(3).
All we need to know is that $\phi$ is continuous, which amounts
to proving the following lemma.\\

\begin{lemma}
Let \((\alpha_n)\) be a sequence of irrational rays with common basepoint converging to
another irrational ray $\alpha$.  Then for every \(B\in\Itin\alpha\) we have \(B\in\Itin\alpha_n\)
for large enough $n$.
\end{lemma}

\begin{proof}
Write \(\Itin\alpha=[B_1,B_2,...]\) and choose \(k\ge 1\).
Then \(B_{k+1}\) is a neighborhood of \(\alpha(t)\) for some time
$t$, which means that for large enough $n$ \(\alpha_n(t)\in B_{k+1}\).
Since \(\alpha_n|_{[0,t]}\) begins in $B_1$ and ends in $B_{k+1}$,
\fullref{le:itins}(1) tells us that it must enter $B_k$.
\end{proof}

\begin{corollary}
\label{co:irrationalmap}
The natural map \(\phi: IX\to\partial\N\) determined by itineraries is an irrational map.
\end{corollary}

\section{CAT(0) Knot Groups}
\label{sec:knots}
\subsection{Preliminary Definitions}
Before we begin discussing knot groups we present some standard terminology concerning CAT(0) groups.
Greater detail is given by Bridson and Haefliger \cite[Chapter II.6]{BH}.
Let $h$ be an isometry of a CAT(0) space $Z$.  If there is a geodesic line $L$ such that
$h$ restricts to a non-trivial translation of $L$, then $L$ is called an \textit{axis} of $h$.  For a point \(z\in L\)
the sequence \((h^nz)_{n=1}^\infty\) converges in $Z\cup\partial Z$ to one of the two boundary points of $L$; we call
that boundary point $h^\infty$.  In fact, given any \(z\in Z\) the sequence \((h^nz)\) converges to the same
point \(h^\infty\).\\

The \textit{minset} of an element \(g\in G\), written \(\Min g\), is the subspace of $Z$ where the map \(z\mapsto d(z,gz)\)
achieves its minimum.
If there is an element \(h_0\) in the center of $G$, then \(\Min h_0\) is the union of axes of $h_0$ and
\(\Min h_0\) splits as a CAT(0) product \(Z'\times\R\).  In this structure, the axes of
\(h_0\) are the geodesic lines \(\{z\}\times\R\).\\

There are two notions of angles in CAT(0) spaces.
The first is the \textit{Alexandrov angle}.
Given two geodesics (segments or rays) $\alpha$ and $\beta$ with the same
initial point $p$, the \textit{Alexandrov angle} between them is the angle between their initial
velocities (see \cite[Definition I.1.12]{BH}) and is denoted by
\(\angle_p(\alpha,\beta)\), or \(\angle_p(a,b)\) if $a$ and $b$ are points on $\alpha$ and $\beta$
other than $p$.\\

The other notion of an angle is the \textit{Tits angle}.
Given two points $\eta$ and $\zeta$ in the boundary of a CAT(0) space $Z$ the \textit{Tits angle} or
\textit{Tits distance} between them is defined by
\[
	\angle_{Tits}(\eta,\zeta)
=
	\sup
	\Bigl\{
		\angle_p\bigl(\overrightarrow{p\eta},\overrightarrow{p\zeta}\bigr)
	\Big|
		p\in Z
	\Bigr\}
\]
where \(\overrightarrow{p\eta}\) and \(\overrightarrow{p\zeta}\) denote the geodesic rays
emanating from $p$ going out to $\eta$ and $\zeta$ respectively.
In the Euclidean plane the two notions agree: that is, for any point $p$ and geodesic rays
$\alpha$ and $\beta$ emanating from $p$ we have:
\[
	\angle_p(\alpha,\beta)=\angle_{Tits}\bigl(\alpha(\infty),\beta(\infty)\bigr)
\]
In fact, when either angle is less than $\pi$, this equation holds precisely when the convex hull
of the union of the two rays is a flat sector \cite[Corollary II.9.9]{BH}.\\

Finally, the following terminology will be convenient when talking about CAT(0) spaces which
split as a product \(\Gamma\times\R\) where $\Gamma$ is a tree.
For vertices \(v\in\Gamma\) we refer to the lines \(\{v\}\times\R\) as \textit{vertical lines}.
For a geodesic edge path \(\nu\subset\Gamma\) we refer to the subspace \(\nu\times\R\)
as a \textit{vertical strip}.\\

\subsection{Knot Groups of Torus Knots}
\label{sec:torusknots}
A torus knot is a knot which lives in a torus.  Specifically, given a relatively prime pair $(p,q)$ we let
$K=K_{p,q}$ be an imbedding \(S^1\incl T^2\subset S^3\) which wraps the circle $p$ times around one direction
of $T^2$ and $q$ times around the other.  It follows from the Van Kampen theorem that the fundamental group
$G$ of the complement \(S^3-K\) is presented by:
\[
	\bigl<a,b\big|a^p=b^q\bigr>
\]
The center of this group is generated by the element \(a^p=b^q\), which we will denote by $\tau$.
Another important group element is the element which represents a meridianal loop in $S^3$ around $K$,
which we will call $\omega$.  By making appropriate choices, we can get
\[
	\omega = b^na^m\quad\textrm{ where $n,m$ solve the equation }mq+np=1.
\]

\begin{figure}[ht!]
\begin{minipage}[b]{0.45\linewidth}
\centering
\input{kpi1.tex}
\end{minipage}
\begin{minipage}[b]{0.45\linewidth}
\labellist
\small\hair 2pt
\pinlabel {$\scriptstyle\times 2$} [b] at 140 300
\pinlabel {$\scriptstyle\times 3$} [b] at 565 300
\pinlabel {$S^1\times[0,1]$} [t] at 360 0
\endlabellist
\begin{flushright}
\indexspace
\includegraphics[height=1in,width=2in]{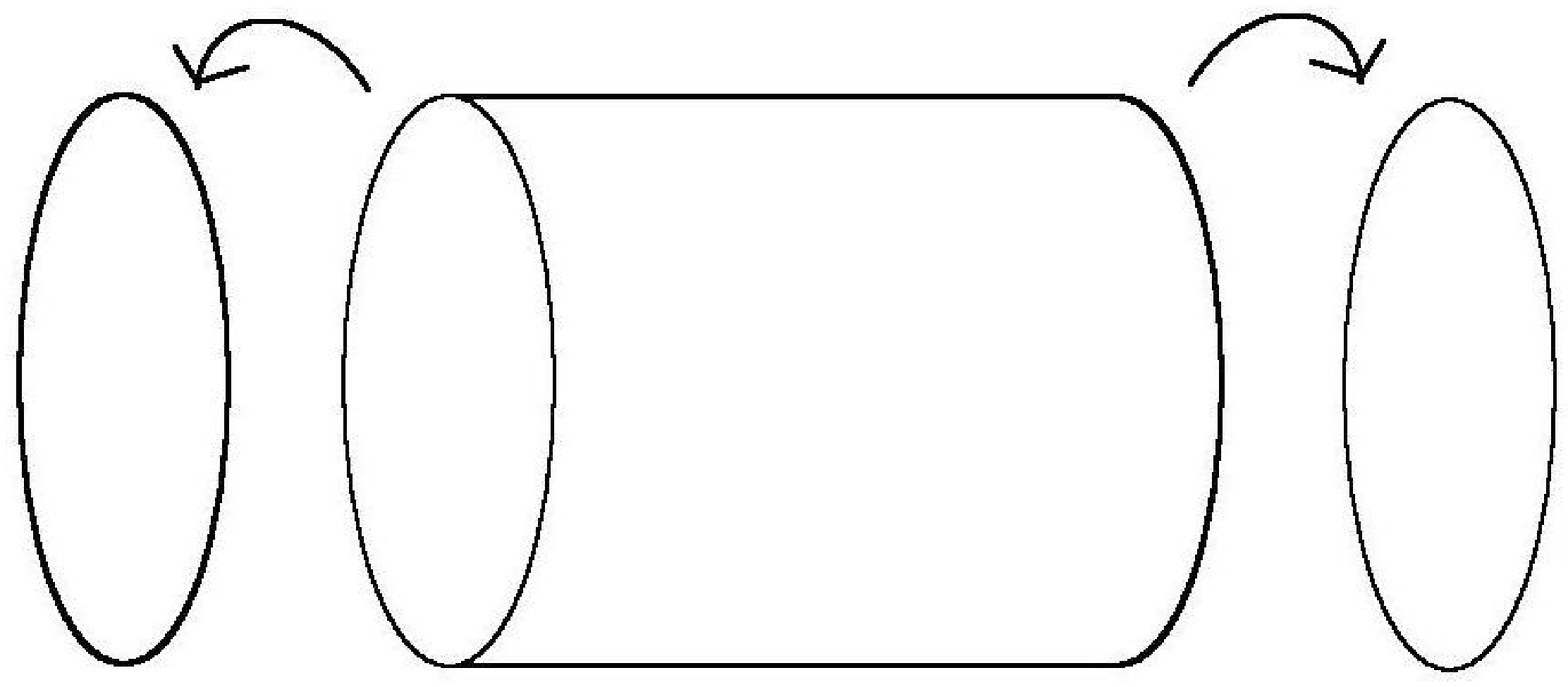}
\end{flushright}
\end{minipage}
\caption{$\overline Y$ for the Trefoil Knot}
\label{fig:Kpi1}
\end{figure}

As in \cite[Example II.11.15(2)]{BH}, we construct a nonpositively curved \(K(G,1)\).
Beginning with a flat rectangle \(\overline{R}=[0,\alpha]\times[0,\beta]\)
of arbitrary dimensions \(\alpha,\beta>0\), we form the quotient space
\(\overline{R}\to\overline{R}/\sim\) where \(\sim\) is generated by the following
three relations:\\
\begin{align*}
	(0,t) &\sim (0,t+\beta/p) \\
	(\alpha,t) &\sim (\alpha,t+\beta/q) \\
	(t,0) &\sim (t,\beta) \\
\end{align*}
This space is nonpositively curved by \cite[Corollary II.11.19]{BH}.  We denote it by \(\overline Y\).
Note that we get the same result if we use the following construction.
Starting with an annulus we glue the two boundary circles to two disjoint circles.
One of the attaching maps wraps the circle $p$
times around itself; the other wraps the circle $q$ times around itself.  
\fullref{fig:Kpi1} shows these two ways to draw \(\overline Y\) for the trefoil knot.
We observe that \(\overline Y\) can be realized topologically as a 2--dimensional spine of the complement of $K$.\\

By the Cartan--Hadamard Theorem (proven by S\,B Alexander and R\,L Bishop in \cite{AB}),
the universal cover \(p: Y\to\overline Y\) is CAT(0).
This CAT(0) $G$--space splits the product \(\Gamma^{p,q}\times\R\)
where \(\Gamma^{p,q}\) denotes the $(p,q)$--biregular tree
\footnote{By ``$(p,q)$--biregular'' we mean the infinite tree whose vertices
alternate in valence between $p$ and $q$.}.  This is the Bass--Serre
tree for the obvious structure as a free product with amalgamation:
\[
	G=\left<a\right>\ast_{\left<\tau\right>}\left<b\right>
\]
\indexspace

\begin{figure}[ht!]
\input{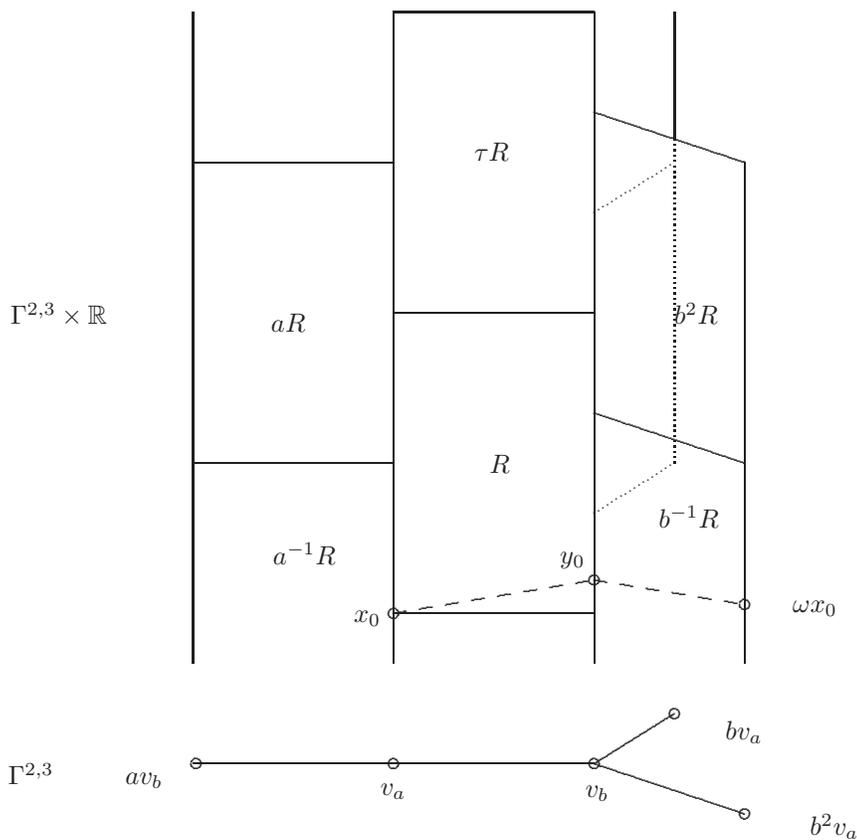}
\caption{$Y$ in the case $p=2$, $q=3$}
\label{fig:Y}
\end{figure}

The action of $G$ on $Y$ is described as follows.  The fundamental chamber is
a lift $R$ of $\overline R$.  The isometry $\tau$ is a vertical translation
by a distance of $\beta$ ($\Min\tau=Y$).
The axis of $a$ is a vertical line containing one side of $R$.
The isometry $a$ is a rotation about this axis followed by a vertical translation
by $\beta/p$.  Similarly, the isometry $b$ is a rotation about its axis followed by
a vertical translation by $\beta/q$.  The action of $G$ on $Y$ is shown in
\fullref{fig:Y}.
In the picture \(v_a\) is the fixed point in \(\Gamma^{p,q}\) of $a$ and
\(v_b\) is the fixed point in \(\Gamma^{p,q}\) of $b$.
We choose as our preferred basepoint a point
\(x_0\) in the axis of $a$.
Also we coordinatize \(Y=\Gamma^{p,q}\times\R\) so that the coordinates of
$x_0$ are \((v_a,0)\) and $\tau$ translates in the positive direction of $\R$.\\

\begin{prop}
\label{prop:hypotenuse}
The geodesic \([x_0,\omega x_0]\) is the hypotenuse of a right triangle
with sides of length $2\alpha$ and $\beta/pq$.
\end{prop}

\begin{proof}
Denote the translation vector of a group element \(g\in G\) in the \(\R\)--coordinate by \(\lambda(g)\)
so that the following hold:
\begin{align*}
	\lambda(\tau) &= \beta\\
	\lambda(a) &= \beta/p\\
	\lambda(b) &= \beta/q\\
\end{align*}
Recalling that $n$ and $m$ satisfy the equation \(mq+np=1\),
we compute:
\begin{align*}
	\lambda(\omega)
&=
	n\lambda(b)+m\lambda(a) \\
&=
	\frac{n\beta}{q}+\frac{m\beta}{p} \\
&=
	\frac{\beta}{pq}
\end{align*}
So the coordinates of $\omega x_0$ in terms of the splitting are:
\[
	\omega x_0
=
	\left(b^{n}v_a,\frac{\beta}{pq}\right)
\]
\end{proof}

Given \(g\in G\) consider the conjugate element \(\omega^g=g\omega g^{-1}\).  Then
\(\Min\omega^g\) is a Euclidean plane of the form \(L\times\R\) where $L$ is the (unique) axis of $\omega^g$ in \(\Gamma^{p,q}\).
We call these Euclidean planes \textit{walls} because they will play the same role as the walls
described in \fullref{sec:origCK}.  The geodesic lines $L$ will be called \textit{wall shadows}.\\

\begin{lemma}
The number of wall shadows containing a given vertex is equal to the valence of the vertex.
The number of wall shadows containing a given edge is 2.\\
\end{lemma}

\begin{proof}
We begin by proving that the number of wall shadows containing a $p$--valent vertex is $p$.
Since $G$ acts transitively on the collection of $p$--valent vertices it suffices to prove this for $v_a$.
Let $L$ be a wall shadow containing $v_a$, say the axis of $\omega^h$.  Then \(L=hL_0\) where
$L_0$ is the axis of $\omega$.  Since \(h^{-1}v_a\in L_0\) there is a $k$ such that
\(\omega^kh^{-1}v_a=v_a\); that is, \(\omega^kh^{-1}\) fixes $v_a$ and is therefore
a power of $a$, say \(\omega^kh^{-1}=a^{-i}\).  Then \(h=a^i\omega^k\) and
\[
	\omega^h=\omega^{a^i}.
\]
There are exactly $p$ conjugate elements of this form giving us $p$ wall shadows
containing the vertex $v_a$.  A similar argument works for $q$--valent vertices.\\

We now prove the second statement of the lemma.  Since $G$ acts transitively on the collection
of edges, every edge is contained in the same number of wall shadows.  Call this number
$\epsilon$.  Consider the number $N$ of pairs \((e,L)\) where $e$ is an edge
of the star of $v_a$ which is contained in the wall shadow $L$.
On one hand since every wall shadow hits two edges of the star of $v_a$, we have \(N=2p\).
On the other hand since every edge is contained in $\epsilon$ wall shadows, we have
\(N=\epsilon p\).  So \(\epsilon=2\).
\end{proof}

\begin{lemma}
The intersection of two wall shadows is at most two edges.
In fact, two wall shadows can contain more than one edge only when $p$ or $q$ is 2.
\end{lemma}

\begin{proof}
It suffices to prove the following claim.

\begin{claim}
Let $v$ be any vertex and $L$ and $L'$ be two wall shadows containing $v$.
Then $L$ and $L'$ share two edges of the star of $v$ iff the valence of $v$ is 2.
\end{claim}

By translating the picture we may assume that $v=v_b$ or $v_a$.  Since the argument is the same
either way, we will assume \(v=v_b\).  Certainly if $q=2$ then $L$ and $L'$ have to share
both edges in the star of $v_b$.  For the converse, assume \(L\cap L'\) contains two edges in the
star of $v_b$.  Without loss of generality assume one of the edges is \([v_a,v_b]\).
Then one of the wall shadows is the axis of $\omega$ and the other is the axis of \(\omega^{b^{-n}}=a^mb^n\).
Say that $L$ is the former and $L'$ is the latter.
The two vertices in the link of $v_b$ hit
by $L$ are $v_a$ and $b^nv_a$.  The two vertices in the link of $v_b$ hit by $L'$ are $v_a$ and
$b^{-n}v_a$.  So $L$ and $L'$ share two edges in the star of $v_b$ only if \(b^nv_a=b^{-n}v_a\)
which happens precisely when \(q|2n\).
Since $q$ and $n$ are relatively prime this is the same as saying that \(q=2\).
This proves the claim.
\end{proof}

Roughly speaking the above lemma tells us that wall shadows bifurcate at odd--valent vertices.
Translated into the language of walls this means several things.\\

\begin{fact}
If $e$ is an edge of \(\Gamma^{p,q}\), then the vertical strip \(e\times\R\) is contained in exactly
two walls.
\end{fact}

\begin{fact}
If $v$ is a $p$--valent [$q$--valent] vertex of $\Gamma^{p,q}$, then the vertical line
\(v\times\R\) is contained in exactly $p$ [$q$] walls.
\end{fact}

\begin{fact}
\label{fact:wallstrip}
Let $W$ and $W'$ be walls.  Then \(W\cap W'\) is either empty, a vertical line or a vertical strip.\\
\end{fact}

Let $\gamma_0$ denote the axis of $\omega$ containing the point $x_0$ and $\overline\gamma$ denote its image in $\overline Y$.
This is a local geodesic loop in $\overline Y$ representing $\omega$.
The $G$--translates of $\gamma_0$ are called \textit{joint lines}, for reasons which will become apparent in the next section.\\

\begin{prop}
Joint lines do not intersect.
\end{prop}

\begin{proof}
Suppose two joint lines intersect at a point $z$.  Without loss of generality we may assume that
one of the joint lines is $\gamma_0$ and that \(z\in R\).  Call the other joint line $\gamma$.
Now if $z$ is in the axis of $a$, then \(\gamma=a^k\gamma_0\) for some $k$.  But this means that \(z=a^kz\),
which is impossible.  A similar argument shows that
$z$ cannot be in the axis of $b$.  Therefore $z$ is in the open vertical strip \((v_a,v_b)\times\R\)
and \(\gamma=\tau^kb^{-n}\gamma_0\) for some $k$.
If $y_0$ is the point at which $\gamma_0$ hits the axis of $b$, then
\[
	z=[x_0,y_0]\cap[x_0',y_0']
\]
where \(y_0'=\tau^kb^{-n}y_0\) and \(x_0'=\tau^ka^mx_0\).
Let \(r: Y\to\R\) denote the projection onto the $\R$--coordinate so that
\(r(x_0)=0\) and \(r(y_0)=\beta/2pq\).
A computation gives \(r(x_0')=r(y_0')+\beta/2pq\).
Therefore:
\[
	0<r(x_0')\le\frac{\beta}{pq}
\]
But \(r(x_0')=i\beta/q\) for some integer $i$, which gives us a contradiction.
\end{proof}

Since the axes of a group element are contained in that element's minset,
the following is true.\\

\begin{fact}
Two joint lines are parallel iff they are contained in the same wall.
\end{fact}

Here we are using the word ``parallel'' in the strong Euclidean sense; that is, when we say two lines
are parallel we mean that their convex hull is a flat strip.  The proof of this next proposition is
an immediate consequence of \fullref{prop:hypotenuse}.\\

\begin{prop}
\label{prop:theta}
We can choose the dimensions $\alpha$ and $\beta$ of $\overline R$ so that
$\overline\gamma$ has length 1 and
\[
	\angle_{Tits}(\omega^\infty,\tau^\infty)
=
	\angle_{x_0}(\omega x_0,\tau x_0)
=
	\theta
\]
for any \(0<\theta<\pi/2\) we choose.
This done, joint lines form angle $\theta$ with vertical lines.
\end{prop}

We close this section with two propositions.  The first is recorded to demonstrate the need for
Wilson's work \cite{Wi} as noted at the end of \fullref{sec:origCK}.  The second shows that
knot groups of Torus knots are rigid, a fact noted in the introduction.\\

\begin{prop}
\label{prop:nopiover2}
It is impossible to construct a CAT(0) $G$--space $Y$ in such a way that
\[
	\angle_{Tits}(\omega^\infty,\tau^\infty)=\pi/2.
\]
\end{prop}

\begin{proof}
Suppose $Y$ is a CAT(0) $G$--space.  Without loss of generality
we may assume \(Y=\Min\tau\) and splits as \(Y=\Gamma'\times\R\).
Then we can define $\lambda$ as in \fullref{prop:hypotenuse} and compute:
\begin{align*}
	\lambda(\omega)
&=
	n\lambda(b) + m\lambda(a) \\
&=
	\frac{\lambda(\tau)}{pq} \\
&\neq
	0
\end{align*}
Choosing \(x_0\in\Min\omega\), we have:
\[
	\angle_{Tits}(\omega^\infty,\tau^\infty)
=
	\angle_{x_0}(\omega x_0,\tau x_0) \\
<
	\frac\pi2
\]
\end{proof}

\begin{prop}
\label{prop:rigidtorusknots}
$G$ is rigid.
\end{prop}

\begin{proof}
By a result of Bedenikovic, Delgado and Timm \cite[Lemma 4.2]{BDT},
we know that $\overline Y$ has a nontrivial self-cover.
By \cite[Theorem 5.2]{BDT}, it has a finite cover \(S^1\times\mathcal{G}\to\overline Y\)
where $\mathcal{G}$ is a finite graph.  Therefore $G$ contains the group \(F\times\Z\)
as a finite index subgroup for some finitely generated free group $F$.\\

Thus if we are given any CAT(0) $G$--space $Y$, the induced action of \(F\times\Z\) on $Y$ as a subgroup 
is cocompact and hence geometric.  Therefore any CAT(0) boundary of $G$ is also a boundary of \(F\times\Z\).
Applying the result of Bowers and Ruane \cite{BR}, we get that every boundary of $G$ is homeomorphic to the suspension of
a cantor set.
\end{proof}

\subsection{Knot Groups of Connected Sums of Torus Knots}
\label{sec:connectedsumsoftorusknots}
Take two relatively prime pairs \((p_\pm,q_\pm)\) and form the corresponding torus knots
\(K_\pm\subset S^3\).  Denote the fundamental group of the complement of
\(K_\pm\) by \(G_\pm\) and let \(\omega_\pm\in G_\pm\) denote the group element representing a
meridianal loop as in \fullref{sec:torusknots}.
Let $K$ be a connected sum \(K_-\#K_+\) and set
\[
	G=\pi_1(S^3-K)=G_-\ast_\Z G_+
\]
where \(\Z\incl G_\pm\) is given by \(1\mapsto\omega_\pm\).
Fixing \(\theta\in(0,\pi/2)\), form the $K(G_\pm,1)$ prescribed in the previous
section and call it \(\overline Y_\pm\).
Construct it so that the local geodesic \(\overline{\gamma_\pm}\subset\overline Y_\pm\) corresponding to
the group element \(\omega_\pm\) has length 1, and in the universal covers
\(Y_\pm\) of \(\overline Y_\pm\) joint lines form angle $\theta$ with vertical lines.
Glue \(\overline Y_-\) to \(\overline Y_+\) along an isometry \(\overline\gamma_-\cong\overline\gamma_+\) to form a
nonpositively curved $K(G,1)$ which we call \(\overline X\).
Let \(p: X\to\overline X\) be the universal covering projection.
Then $X$ is a CAT(0) $G$--space.
Since \(G_\pm\) both inject into $G$, the path components of \(p^{-1}(\overline Y_\pm)\)
are isometric copies of \(Y_\pm\).
We call these path components \textit{natural blocks}.
It is easy to see that the collection of natural blocks gives us a block structure on $X$
in the sense of \fullref{sec:NervesandItineraries}.  Thus the nerve $\N$ of the collection
of natural blocks is a tree and we may talk about the itinerary between two natural blocks or
the itinerary of a geodesic.
We call an itinerary in terms of natural blocks a \textit{natural itinerary} and use the notation
\(\Itin_\N\).\\

Now this ``natural block structure'' is different from the block structure of Croke and Kleiner's
construction in \cite{CK}.
Here natural blocks do not intersect at walls (Euclidean planes) but at joint lines.
We will see, however, that the boundary of our construction has the same essential structure as the boundary
in \cite{CK}.  To prove this we will need to introduce another type of block.\\

\begin{definition}
Given a joint line $\gamma$ we define the \textit{joint block of $\gamma$} to be the convex
hull of all joint lines $X$ which are parallel to $\gamma$.  This done, we define a ``new nerve'' $\widehat\N$
with the following properties:\\
\begin{enumerate}
\item Vertices \(\widehat v_B\) correspond to blocks $B$ of $X$ (joint and natural).\\
\item An edge \([\widehat v_{B_1},\widehat v_{B_2}]\) is included whenever \(B_1\cap B_2\) is a wall.\\
\end{enumerate}
When (2) holds we will say that $B_1$ \textit{neighbors} $B_2$.\\
\end{definition}

A word of warning:  When we call $\widehat\N$ a ``nerve'' we do not mean it in the
same sense as used in \fullref{sec:NervesandItineraries}.  We mean here that it is the correct
analogue of the previous notion of a nerve in this context.  In \cite{CK} at most two blocks
could intersect simultaneously and then their intersection was precisely a wall.  Here there are
many intersections which are not being recorded; for example, every point is in at least three blocks,
possibly more.  This fact will cause some difficulty for us in \fullref{sec:itinsinwidehatN} when we need
to redefine itineraries of terms of $\widehat\N$.\\

\section{The Main Theorem}

In order to apply the strategies of Croke and Kleiner \cite{CK} and Wilson \cite{Wi}, will need to prove that if we take
the collection of all blocks, both joint and natural, together with this ``new nerve'' $\widehat\N$, then
Theorems A--E from \fullref{sec:origCK} remain valid.  Restated in this context the theorems will
be labeled A$'$--E$'$.\\

\subsection{Joint Blocks}
If \([B_1,...,B_n]\) is a natural itinerary, then we call the list of joint lines
\begin{align*}
	\gamma_1 &= B_1\cap B_2\\
	\gamma_2 &= B_2\cap B_3\\
		&\vdots \\
	\gamma_{n-1} &= B_{n-1}\cap B_n
\end{align*}
the list of joint lines \textit{between \(B_1\) and \(B_n\)}.  If \(\gamma\) and \(\gamma'\) are two joint
lines, then it is easy to see that every geodesic which begins on \(\gamma\) and ends on \(\gamma'\) has the same itinerary.
If that itinerary is \([B_1,...,B_n]\), and \(\gamma_1,...,\gamma_{n-1}\) is the list of joint lines between
$B_1$ and $B_n$, then we also call \(\gamma_1,...,\gamma_{n-1}\) the
the list of joint lines \textit{between $\gamma$ and $\gamma'$}.\\

\begin{lemma}
(The Joint Line Lemma)
Let \(\gamma\) and \(\gamma'\) be parallel joint lines.
Then every joint line between $\gamma$ and $\gamma'$ is also parallel to $\gamma$ and $\gamma'$.
\end{lemma}

\begin{proof}
Parameterize \(\gamma,\gamma':\R\to X\) to have unit speed and let $\gamma_0$ be a joint line
between $\gamma$ and $\gamma'$.  Then every geodesic which begins on $\gamma$ and ends on $\gamma'$ must pass through
$\gamma_0$ (\fullref{le:itins}(2)).
In particular, for \(k\in\Z\) the geodesic
\([\gamma(k),\gamma'(k)]\) intersects $\gamma_0$ at some point $z_k$.  Since \(\gamma\|\gamma'\),
\(d(\gamma(k),\gamma'(k))\) is constant and \(z_k\) remains asymptotic to $\gamma$ and $\gamma'$
as \(k\to\pm\infty\).  It follows that $\gamma_0$ is indeed parallel to these.
\end{proof}

\begin{prop}
Let $B_J$ be the joint block of a joint line $\gamma_0$.  Then:\\
\begin{enumerate}
\item \(B_J\cong\Gamma^4\times\R\) where \(\Gamma^4\) is the 4--valent tree.\\
\item The joint lines parallel to \(\gamma_0\) are precisely the vertical lines in \(B_J\).\\
\item If a joint line \(\gamma\subset B_J\), then \(\gamma\|\gamma_0\).\\
\end{enumerate}
\end{prop}

\begin{proof}
Let \(D^0=\gamma_0\), and for each \(n>0\), let \(D^n\) denote the union of \(D^{n-1}\)
along with all walls intersecting \(D^{n-1}\) at a joint line.  In addition, we define
\(D^\infty=\bigcup_{i=0}^\infty D^i\).  Since two nonintersecting lines in a common
Euclidean plane must be parallel, we see that \(D^n\) splits as \(T^n\times\R\) where
each $T^n$ is a tree, constructed as follows: We begin with $T^0$, which is just a point.
$T^1$ is the union of two lines glued together at a single point $z_0$.  To form $T^2$, we
glue four new lines to $T^1$ at four points \(z_1,...,z_4\) in the four components of \(T^1-z_0\).
To form $T^3$, we glue twelve new lines to $T^2$ at twelve points in the twelve unbounded components
of \(T^2-\{z_1,...,z_4\}\), and so on.  The limit $\Gamma^4$ of this increasing sequence of trees is
an infinite 4--valent tree.  Thus we get \(D^\infty=\Gamma^4\times\R\).
Furthermore, we see from the construction that the joint lines in \(D^\infty\) are
precisely the vertical lines, and that all of these are parallel to $\gamma_0$.
So the proposition will follow if we show that \(D^\infty=B_J\).\\

Certainly \(D^\infty\) is the convex hull of the collection of joint lines parallel
to \(\gamma_0\) which are contained in \(D^\infty\).  What we need to know is
that \textit{all} joint lines parallel to $\gamma_0$ are contained in \(D^\infty\).
We prove this here:
Let \(\gamma\) be a joint line parallel to \(\gamma_0\), and \(\gamma_1,...,\gamma_{n-1}\)
be the list of joint lines between $\gamma_0$ and $\gamma$.  Since \(\gamma_1\|\gamma_0\) and
these two are in a common natural block, it follows
that they are in a common wall and that \(\gamma_1\subset D^1\).
In general, since \(\gamma_i\|\gamma_{i+1}\) and these
two joint lines are in a common natural block, they are in a common wall and therefore
\(\gamma_{i+1}\subset D^{i+1}\).  So \(\gamma_{n-1}\subset D^{n-1}\) and \(\gamma\subset D^n\).
\end{proof}

\begin{rem}
This proposition corresponds to the group theoretic fact that the stabilizer of a joint
block is \([\Z\ast\Z]\times\Z\).  For example, if \(\gamma_0\) is the joint line containing \(x_0\)
then the stabilizer of \(B_J\) is \([\left<\tau_-\right>\ast\left<\tau_+\right>]\times\left<\omega\right>\)
where $\omega$ translates \(B_J\) in the $\R$--direction and \(\left<\tau_-\right>\ast\left<\tau_+\right>\)
acts geometrically on \(\Gamma^4\).  Here $<\tau_\pm>$ denote the centers of \(G_\pm\).
\end{rem}

Two distinct blocks neighbor each other iff one is joint, the other is natural, and the two share a joint line.
For a joint block \(B_J\) let \(\C(B_J)\) denote the collection of natural blocks which neighbor \(B_J\)
and \(\N(B_J)\) denote the full subgraph of $\N$ spanned by the vertices \(\{v_{B_N}|B_N\in \C(B_J)\}\).\\

\begin{lemma}
\label{le:CandNofBJ}
Let $B_J$ and $B'_J$ be distinct joint blocks.  Then:\\
\begin{enumerate}
\item If \([v_{B_N},v_{B'_N}]\) is an edge of \(\N(B_J)\), then the joint line \(B_N\cap B'_N\) is in \(B_J\).\\
\item If \(B_N,B'_N\in \C(B_J)\), then \(\Itin_\N[B_N,B'_N]\subset \C(B_J)\).\\
\item \(\bigl|\C(B_J)\cap \C(B'_J)\bigr|\le 1\).\\
\end{enumerate}
\end{lemma}

\begin{proof}
(1) Let $\gamma$ and $\gamma'$ be joint lines of $B_N$ and $B'_N$ which are in $B_J$.  Since \(B_N\cap B'_N=\gamma_0\)
is a joint line, $\gamma_0$ is the only joint line between $\gamma$ and $\gamma'$.
It follows from the joint line lemma that $\gamma_0$ is parallel to $\gamma$ and $\gamma'$ and must therefore
also be in $B_J$.\\

(2) Again, let $\gamma$ and $\gamma'$ be joint lines of $B_N$ and $B'_N$ which are in $B_J$, and write
\(\Itin_\N[B_N,B'_N]=[B_1,...,B_k]\).  Then for \(1\le i<k\), the joint line lemma tells us that the joint
lines \(B_i\cap B_{i+1}\) are all in \(B_J\).  So for \(1\le i\le k\), every $B_i$ shares a joint line with
\(B_J\).\\

(3) Suppose \(|\C(B_J)\cap \C(B'_J)|>1\).
Since \(\N(B_J)\) and \(\N(B'_J)\) are convex, \(\N(B_J)\cap \N(B'_J)\) must contain an edge \([v_{B_N},v_{B'_N}]\);
by (1), the joint line \(B_N\cap B'_N\) is in both \(B_J\) and \(B'_J\), which is a contradiction.
\end{proof}

\subsection{\texorpdfstring{Itineraries in $\widehat\N$}{Itineraries in the ``New Nerve''}}
\label{sec:itinsinwidehatN}

Our goal here is to show that $\widehat\N$ is a tree (Theorem A$'$)
and to define itineraries in terms of \(\widehat\N\).\\

\begin{lemma}
\label{le:itinerarychain}
Let \([\widehat v_{B_1},...,\widehat v_{B_n}]\) be an edge path in $\widehat\N$ with no backtracking
such that $B_1$ and $B_n$ are natural blocks.
Then
\begin{align*}
	\Itin_\N[B_1,B_n]
&=
	\Itin_\N[B_1,B_3]\cup\Itin_\N[B_3,B_5]\cup \\
&
	\;.\;.\;.\; \cup\Itin_\N[B_{n-4},B_{n-2}]\cup\Itin_\N[B_{n-2},B_n]
\end{align*}
where
\[
	\Itin_\N[B_{i-1},B_{i+1}] \subset\C(B_i)
\]
for even \(1<i<n\).
\end{lemma}

\begin{proof}
We prove this by induction on $n$.  When \(n=1\), there is nothing to show, and when \(n=3\), we simply
note that \(\N(B_2)\) is convex.  Assume \(n\ge 5\),
and let $\eta$ denote the geodesic edge path in $\N$ from \(v_{B_1}\) to \(v_{B_{n-2}}\);
since \(B_{n-4}\neq B_{n-2}\), the last edge of $\eta$ is in \(\N(B_{n-3})\) (by induction).
It follows that \(\eta\cap \N(B_{n-1})=\{v_{B_{n-2}}\}\).  Hence, if \(\eta'\) is
the geodesic edge path from \(v_{B_{n-2}}\) to \(v_{B_n}\), then since \(\eta'\subset \N(B_{n-1})\),
the edge path \(\eta\cup\eta'\) has no backtracking and must be the geodesic edge path in $\N$ between
\(v_{B_1}\) and \(v_{B_n}\).
\end{proof}

\begin{ThmA'}
\(\widehat\N\) is a tree.
\end{ThmA'}

\begin{proof}
It follows from \fullref{le:CandNofBJ}(3) that $\widehat N$ has no squares.
Thus, any non-nullhomotopic loop in \(\widehat N\) must have length at least 6.
Suppose \([\widehat v_{B_1},...,\widehat v_{B_n}]\) is such a loop with no backtracking
where \(B_1=B_n\) is natural.  Then by the previous lemma, the first edge in the geodesic
edge path in $\N$ from $v_{B_1}$ to $v_{B_3}$ is in both \(\N(B_2)\) and \(\N(B_{n-1})\),
giving us a contradiction.
\end{proof}

\begin{prop}
\label{prop:hitswalls}
Suppose \([\widehat v_{B_1},...,\widehat v_{B_n}]\) is a geodesic edge path in \(\widehat\N\)
and $\alpha$ is a geodesic segment which begins in \(B_1\) and ends in \(B_n\).  Then
then $\alpha$ is covered by the collection of blocks \(\{B_k\}_{k=1}^n\) and passes through
every block $B_k$ and wall $B_k\cap B_{k+1}$.
\end{prop}

\begin{proof}
First of all, assume $B_1$ and $B_n$ are both natural blocks and use \fullref{le:itinerarychain} to write:
\begin{align*}
	\Itin_\N\alpha
&\subset
	\Itin_\N[B_1,B_n] \\
&=
	\Itin_\N[B_1,B_3]
		\cup
			...
		\cup
	\Itin_\N[B_{n-2},B_n]
\end{align*}
So for every odd \(1<k<n\), $\alpha$ passes through the block \(B_k\).
Let \(\gamma_{k-2}\) denote the joint line at which $\alpha$ leaves the natural block \(B_{k-2}\)
and \(\gamma_{k-1}\) denote the joint line at which $\alpha$ enters the natural block \(B_k\).
The fact that
\(\Itin_\N[B_{k-2},B_k]\subset C(B_{k-1})\) tells us that \(\gamma_{k-2},\gamma_{k-1}\subset B_{k-1}\).
For every \(1\le k<n\), let $t_k$ be the time such that \(\alpha(t_k)\in\gamma_k\).
Since \(\gamma_k\subset B_k\cap B_{k+1}\), we see that $\alpha$ hits every such wall.
Furthermore,
\[
	\gamma\bigl([t_k,t_{k+1}]\bigr)\subset B_k
\]
because $B_k$ is convex.  This shows that \(\alpha\subset\bigcup_{k=1}^nB_k\).\\

Now consider the more general case.
If $B_1$ is joint and $B_n$ is natural, choose a natural block \(B_0\) containing the initial point of
\(\alpha\).  This time \fullref{le:itinerarychain} gives us:
\begin{align*}
	\Itin_\N\alpha
&\subset
	\Itin_\N[B_0,B_n] \\
&=
	\Itin_\N[B_0,B_2]
		\cup
	\Itin_\N[B_2,B_4]
			...
		\cup
	\Itin_\N[B_{n-2},B_n]
\end{align*}
As before, since $\alpha$ enters \(B_2\) at a joint line of \(B_1\), we get that $\alpha$ passes through
the wall \(B_1\cap B_2\) and because blocks are convex, we get that \(\alpha\subset\bigcup_{k=1}^nB_k\).
Similar arguments work if $B_1$ is natural and $B_n$ is joint, or if both $B_1$ and $B_n$ are joint.
\end{proof}

We now know that given a geodesic segment (or ray) $\alpha$ there is a (possibly infinite) geodesic edge path
\([\widehat v_{B_1},...,\widehat v_{B_n}]\) such that \(\alpha\subset\bigcup_{k=1}^nB_k\).  We define
the $\widehat\N$--itinerary of $\alpha$ to be the list \([B_1,...,B_n]\) where
\([\widehat v_{B_1},...,\widehat v_{B_n}]\) is the shortest such edge path.
We may also write \(\Itin_{\widehat\N}[B_0',B_1']=[B_1,...,B_n]\)
when \([\widehat v_{B_1},...,\widehat v_{B_n}]\) is the
geodesic edge path in $\widehat\N$ from \(\widehat v_{B_0'}\) to
\(\widehat v_{B_1'}\).\\

There is some danger of confusion here since every geodesic in $X$ has two itineraries: one in terms of $\N$ and
the other in terms of $\widehat\N$.  We already have a notion of rational and irrational rays in terms of
$\N$--itineraries.  We denote the set of points which are ``rational with respect to $\N$'' by \(R_\N X\) and
the set of points which are ``irrational with respect to $\N$'' by \(I_\N X\).
However, in this section we will call a geodesic ray \textit{rational} if its \(\widehat\N\)--itinerary
is finite and \textit{irrational} if its \(\widehat\N\)--itinerary is infinite.
The endpoint of a rational ray is called a \textit{rational point} and
the endpoint of an irrational ray is called an \textit{irrational point}.
To emphasize that by ``rational'' and ``irrational'' we mean in terms of $\widehat\N$, we will denote
the set of rational points of \(\partial X\) by \(R_{\widehat\N}X\) and
the set of irrational points of \(\partial X\) by \(I_{\widehat\N}X\).\\

\begin{prop}
\label{prop:determineditinerary}
If $\alpha$ and $\beta$ are two irrational geodesic rays whose $\widehat\N$--itineraries
eventually coincide, then their $\N$--itineraries are also infinite and also eventually coincide.
\end{prop}

\begin{proof}
Suppose we have two irrational geodesic rays $\alpha$ and $\beta$ whose itineraries eventually coincide.
Write
\[
	\Itin_{\widehat\N}\alpha=[B_1,B_2,...]
\]
and:
\[
	\Itin_{\widehat\N}\beta=[B'_1,B'_2,...]
\]
Then there are $m,n>2$ such that
\(B_{m+i}=B'_{n+i}\) for \(i\ge 0\).  Choose $m$ and $n$ so that \(B_m=B'_n\) is a natural block.
The fact that \(m,n>2\) guarantees that $\alpha$ and $\beta$ do not begin in this block; hence
\(B_m\in\Itin_\N\alpha\cap\Itin_\N\beta\).  In fact, \(B_{m+2i}\in\Itin_\N\alpha\cap\Itin_\N\beta\)
for every \(i\ge 0\).  Therefore \(\Itin_\N\alpha\cap\Itin_\N\beta\) contains the infinite
sequence of blocks \(\{B_m,B_{m+1},B_{m+2},...\}\).
\end{proof}

\begin{corollary}
\hspace{0cm}
\label{co:IRrelationships}
\begin{enumerate}
\item \(I_{\widehat\N}X\subset I_\N X\)\\
\item \(R_\N X\subset R_{\widehat\N}X\)\\
\end{enumerate}
\end{corollary}

\begin{rem}
The above inclusions are strict.
A geodesic ray which stays in the same joint block but does not
stay in any wall will have an infinite $\N$--itinerary but finite $\widehat\N$--itinerary.
\end{rem}

\subsection{\texorpdfstring{The Boundary of $X$}{The Boundary of X}}
\label{sec:KnotGroupBoundary}
We do not yet know that $R_{\widehat\N}X$ is precisely the union of block boundaries.
For this we need to know \fullref{le:irrationaldistance} in the new context.
The proof is the same except that we replace the claim with the following lemma.\\

\begin{lemma}
There is a \(\delta>0\) such that for natural blocks $B$ and $B'$,
if \(d_{\widehat\N}(\widehat v_{B},\widehat v_{B'})\ge 4k\) then \(d(B,B')\ge k\delta\).
\end{lemma}

\begin{proof}
Let $\delta$ be the minimum positive distance between joint lines in $X$.
Let $B$ and $B'$ be natural blocks, \(\Itin_\N[B,B']=[B_0,...,B_n]\)
where \(n\ge 4k\), \(x\in B\) and \(x'\in B'\).
Then the geodesic $[x,x']$ passes through
every block \(B_{4i}\) for \(0\le i\le k\) at some point $z_i$.
Furthermore, for \(0\le i<k\) the geodesic \([z_i,z_{i+1}]\) enters the
block \(B_{4i+2}\) at a joint line of the
joint block \(B_{4i+1}\) and leaves at a joint line of the block \(B_{4i+3}\).
Thus we have
\[
	d(x,x')
\ge
	\sum_{i=0}^{k-1}d(z_i,z_{i+1}) \\
\ge
	k\delta.
\]
\end{proof}

\begin{corollary}
\label{co:newRXIX}
$R_{\widehat\N}X$ is the union of block boundaries and $I_{\widehat\N}X$ is its complement.
\end{corollary}

Since every block $B$ splits as \(\Gamma\times\R\) for some tree $\Gamma$
(one of $\Gamma^4$, $\Gamma^{p_-,q_-}$ or $\Gamma^{p_+,q_+}$),
it follows that \(\partial B\) is the suspension of a cantor set.
As mentioned in \fullref{sec:origCK}, the suspension
points are called \textit{poles} and the set of poles is denoted $PB$.\\

\begin{ThmB'}
Let \(B_0\) and \(B_1\) be blocks and $D$ be the distance between the corresponding vertices in
$\widehat N$.  Then:\\
\begin{enumerate}
\item If \(D=1\), then \(\partial B_0\cap\partial B_1=\partial W\) where $W$ is the wall \(B_0\cap B_1\).\\
\item If \(D=2\), then \(\partial B_0\cap\partial B_1=PB_{1/2}\) where \(B_{1/2}\) neighbors both
\(B_0\) and \(B_1\).\\
\item If \(D>2\), then \(\partial B_0\cap\partial B_1=\emptyset\).\\
\end{enumerate}
\end{ThmB'}

\begin{proof}
(1) If \(D=1\), then \(B_0\cap B_1\) is a wall $W$.  That \(\partial W\subset\partial B_0\cap\partial B_1\)
is obvious.  The reverse inclusion follows by the same sort of argument as was used in the joint line lemma:
If \(\alpha_0\subset B_0\) and \(\alpha_1\subset B_1\) are asymptotic geodesic rays, then every
geodesic from $\alpha_0$ to $\alpha_1$ intersects the wall $W$.
Thus we can get a sequence of points in $W$ which remain asymptotic to $\alpha_0$ and $\alpha_1$.\\

(2) If \(D=2\), then there is one vertex between \(\widehat v_{B_0}\) and \(\widehat v_{B_1}\); call it
\(\widehat v_{B_{1/2}}\).
We will show that
\[
	PB_{1/2}
		\subset
	\partial B_0\cap\partial B_1
		\subset
	\partial W_0\cap\partial W_1
		\subset
	PB_{1/2}
\]
where \(W_i=B_{1/2}\cap B_i\) for \(i=0,1\).
The first inclusion is just the fact that vertical lines of \(B_{1/2}\) can be found in \(W_0\) and \(W_1\),
and the second is the same argument as in (1).  For the third inclusion, suppose
\(\alpha_0\subset W_0\) and \(\alpha_1\subset W_1\) are asymptotic geodesic rays and let \(\overline\alpha_0\)
and \(\overline\alpha_1\) be the projections of \(\alpha_0\) and \(\alpha_1\) onto the $\Gamma$--coordinate of
\(B_{1/2}=\Gamma\times\R\).  If \(\overline\alpha_0\) and \(\overline\alpha_1\) are not constant, then since
they are asymptotic they must have infinitely many vertices of $\Gamma$ in common.  In this case \(W_0\cap W_1\)
shares a half-plane, contradicting the fact that \(W_0\cap W_1\) is at most a vertical strip.
So \(\overline\alpha_0\) and \(\overline\alpha_1\) are constant and \(\alpha_0\) and \(\alpha_1\) go to a pole of
\(B_{1/2}\).\\

Finally we show (3) by contradiction: Suppose \(\zeta\in\partial B_0\cap\partial B_1\)
and write \(\Itin_{\widehat\N}[B_0,B_1]=[\overline B_1,...,\overline B_n]\) where \(n=D+1\) by
hypothesis.
By the same argument as in (1) we actually have that
\(\zeta\in\partial B_i\) for every \(1\le i\le n\).  By (2) it follows that
\(\zeta\in PB_i\) for every \(1<i<n\).  But \(PB_2\cap PB_3=\emptyset\) because \(\angle_{Tits}(PB_2,PB_3)=\theta\),
giving us a contradiction!
\end{proof}

\begin{ThmC'}
Let $B$ be a block and \(\zeta\in\partial B\) not be a pole of any neighboring block.  Then $\zeta$ has
a local path component which stays in $\partial B$.
\end{ThmC'}

\begin{proof}
The proof is the same as the proof of Theorem C with one minor exception.  It could be that
a geodesic ray may exit a wall $W$ via a joint line $\gamma$ of another wall.
But by \fullref{fact:wallstrip}, \(W\cap\gamma\) is compact in this case.  So this exception causes no problems.
\end{proof}

As in \cite{CK} we call \(\zeta\in\partial X\) a \textit{vertex}
if there is a local path component $V$ of $\zeta$ and a local path component $V'$ of an actual pole
$\zeta'$ and a homeomorphism \((V,\zeta)\approx(V',\zeta')\).  A point of \(\partial X\)
is a vertex iff it has a local path component homeomorphic to the open cone on the cantor set via a homeomorphism
which takes $\zeta$ to the cone point.
A path in \(\partial X\) is called \textit{safe} if it passes through
vertices at only finitely many times.  Since $R_{\widehat\N}X$ is just the union of block
boundaries (\fullref{co:newRXIX}), Theorem C$'$ tells us that the only vertices in $R_{\widehat\N}X$
are poles.\\

\begin{ThmD'}
$R_{\widehat\N}X$ is the unique dense safe path component of $\partial X$.
\end{ThmD'}

\begin{proof}
The proof that $R_{\widehat\N}X$ is a safe path component is exactly the same as the
proof of \cite[Lemma 6]{CK}.  The fact that $R_{\widehat\N}X$ is dense follows from
\fullref{le:denseRX} and the fact that \(R_{\widehat\N}X\supset R_{\N}X\).
Now the other safe path components are contained in the path components of $I_{\widehat\N}X$.
Recall that \fullref{co:irrationalmap} provided us with a map \(\phi: I_{\N}X\to\partial\N\)
which is ``irrational with respect to $\N$''.  By \fullref{prop:hitswalls} we know that
the restriction $\widehat\phi$ of $\phi$ to \(I_{\widehat\N}X\) is ``irrational with respect to
$\widehat\N$''.  Since \(\widehat\phi\) takes safe path components to points and no point of \(\im\widehat\phi\)
is dense in \(\im\widehat\phi\) it follows that no safe path component of \(I_{\widehat\N}X\) is dense
in \(I_{\widehat\N}X\).
\end{proof}

\begin{ThmE'}
Let $B$ be a block.\\
\begin{enumerate}
\item The union of boundaries of walls of $B$ is dense in \(\partial B\).\\
\item The closure of the set of poles of neighboring blocks is the same as the set
of points of \(\partial B\) which are a Tits distance of $\theta$ from a pole
of $B$.\\
\end{enumerate}
\end{ThmE'}

\begin{proof}
Let $\alpha$ be any geodesic ray in \(B=\Gamma\times\R\) and \(\overline\alpha\) be its projection onto the
$\Gamma$--coordinate.  Let \(\overline\alpha(t_1),\overline\alpha(t_2),...\) be the sequence of non-bivalent
vertices through which $\overline\alpha$ passes.  Then for every \(n\ge 1\) there is a wall \(W_n\) such
that \(\alpha([0,t_n])\cap W_n=\alpha(t_n)\).  Thus we may bifurcate $\alpha$ at \(\alpha(t_n)\) to get a
(probably new) ray $\alpha_n$ which agrees with $\alpha$ up to time $t_n$ and then stays in $W_n$.
This proves (1).\\

For (2) assume that in the above setup we have \(\angle_{Tits}(\alpha(\infty),\zeta)=\theta\)
where $\zeta$ is a pole of $B$.  This means that $\alpha$ hits vertical lines of $B$ at an angle of $\theta$.
Since $\alpha$ enters the wall $W_n$ at time $t_n$ we have two choices for $\alpha_n$.  For one  
of these choices we will have \(\alpha_n([t_n,\infty))\) 
parallel to the non-vertical lines in $W_n$.  Then \(\alpha_n(\infty)\) will be a pole of a neighboring
block.
\end{proof}

\begin{Thm1}
The knot group $G$ of any connected sum of two non-trivial torus knots has
uncountably many CAT(0) boundaries.
\end{Thm1}

\begin{proof}
We sketch here the key argument of \cite{Wi}.
For \(0<\theta<\pi/2\) construct \(X_\theta\) as above.
Now suppose we have \(0<\theta_1,\theta_2<\pi/2\) and a homeomorphism
\(h: \partial X_{\theta_1}\to\partial X_{\theta_2}\).  Since $h$ takes vertices to vertices
it follows from Theorem D$'$ that it takes $R_{\widehat\N}X_{\theta_1}$ to
$R_{\widehat\N}X_{\theta_2}$.  From here it is not hard to see that $h$ takes poles
to poles, block boundaries to block boundaries and wall boundaries to wall boundaries.
Let $W_1$ be a wall in $X_{\theta_1}$ and $W_2$ be the wall of $X_{\theta_2}$ such
that \(h(\partial W_1)=\partial W_2\).
Using Theorem E$'$ and a proof by induction we find sequence of points
\((z_k)_{k=0}^\infty\subset\partial W_1\) such that
\begin{align*}
	\angle_{Tits}(z_k,z_{k+1}) &= \theta_1 \\
\textrm{and\quad}
	\angle_{Tits}\bigl(h(z_k),h(z_{k+1})\bigr) &= \theta_2. \\
\end{align*}
If $\theta_1$ is a rational multiple of $\pi$, then \(\{z_k\}\) is a finite set
and we can use a counting argument to prove that \(\theta_1=\theta_2\).
If $\theta_1$ is not a rational multiple of $\pi$, then \(\{z_k\}\) is a dense subset of $\partial W_1$
and the same argument no longer works.
Wilson's solution is to use the sequences \((z_k)\) and \((h(z_k))\)
to define two nonstandard orderings of the natural numbers denoted \(\prec_1\) and \(\prec_2\)
such that \(\prec_1\) is equivalent to \(\prec_2\).
She then uses a technical argument to show that this fact implies that
\(\theta_1=\theta_2\).
\end{proof}

\end{document}